# Infimal Convolution and Duality in Convex Optimal Control Problems with Second Order Evolution Differential Inclusions

**Elimhan N. Mahmudov**
*Department of Mathematics, Istanbul Technical University, Istanbul, Turkey,*
*Azerbaijan National Academy of Sciences Institute of Control Systems, Baku, Azerbaijan.*
elimhan22@yahoo.com

**Abstract**. The paper deals with the optimal control problem described by second order evolution differential inclusions; to this end first we use an auxiliary problem with second order discrete and discrete-approximate inclusions. Then applying infimal convolution concept of convex functions, step by step we construct the dual problems for discrete, discrete-approximate and differential inclusions and prove duality results. It seems that the Euler-Lagrange type inclusions are "duality relations" for both prımary and dual problems and that the dual problem for discrete-approximate problem make a bridge between them. Finally, relying to the method described within the framework of the idea of this paper a dual problem can be obtained for any higher order differential inclusions. In this way relying to the described method for computation of the conjugate and support functions of discrete-approximate problems a Pascal triangle with binomial coefficients, can be successfully used for any "higher order" calculations.
 **AMS Subject Classifications**, 34A60, 49N15, 49M25, 90C46.

## 1. Introduction

Many extremal problems, for example, as classical problems of optimal control, differential games, models of economic dynamics, macroeconomic problems, etc. are described in terms of set-valued mappings and form a component part of the modern mathematical theory of controlled dynamical systems and mathematical economics [3-8,10,11,13,15, 26]. The first and second-order ordinary and partial differential inclusions, naturally arising from certain physical and control problems, have attracted the attention of many researchers, and as a result, various qualitative problems, including the existence results have been considered by many authors (see [2,7,12,16-18,21,23,24,28,29] and references therein). In the paper [7]are studied the time optimal control problem with endpoint constraints for a class of differential inclusions that satisfy mild smoothness and controllability assumptions. The paper [23] concerns optimal control of discontinuous differential inclusions of the normal cone type governed by a generalized version of the Moreau sweeping process with control functions acting in both nonconvex moving sets and additive perturbations.

 In the papers [18-22], for optimal control problems of higher order discrete processes and differential inclusions (DFIs) with the use of locally adjoint mappings (LAMs) the necessary and sufficient conditions of optimality are formulated. Along with these the duality theory plays a fundamental role in the analysis of optimization and variational problems. The reader can refer to [1,8,9,25,27] and their references for more details on this topic. It not only provides a powerful theoretical tool in the analysis of these problems, but also paves the way to designing new algorithms for solving them. A key player in any duality framework is the Legendre-Fenchel conjugate transform. Often, duality is associated with convex problems, yet it turns out that duality theory also has a fundamental impact even on the analysis of nonconvex problems. The work [17] is devoted to optimization of so-called first-order partial DFIs in the gradient form on a square domain. In the Euler–Lagrange form, necessary and sufficient conditions are derived for the discrete-approximate and partial DFIs, respectively. The duality theorems are proved and duality relation is established. Use of infimal convolution throughout this work plays a key role in proofs of duality results for problems with second order DSIs and DFIs. The aim in the work [9] is to establish conditions under which strong duality can be guaranteed. To this purpose, even convexity and properness are a compulsory requirement over the involved functions in the primal problem.

In the present work, the optimality conditions for a second order discrete inclusions (DSIs) and DFIs together with their duality approach were considered for the first time. The construction of duality for a second order DFIs is accompanied first by duality of second order DSIs and then discrete-approximate problems, where the duality of the last problem should be expressed in terms of the resulting difference operators. In turn, for the duality problem of a discrete-approximate problem of the second order, skilful computations of adjoint and support functions are required. Consequently, the key to our success is the formulation of the Lemmas 4.1–4.3 and Propositions 4.1–4.3, without which it is hardly ever possible to establish any duality to the problem with second order DSIs and DFIs. To the best of our knowledge, there are a few papers (see [9,16,17,25] and references therein) devoted to duality problems of first order DFIs. Building on these results, we then treat dual results according to the dual operations of addition and infimal convolution of convex functions [1,9,11,14,15].

Thus, the present paper is dedicated to one of the difficult and interesting fields – construction of duality of optimization problems with second order ordinary discrete and DFIs. The posed problems and their dualities are new. The paper is organized in the following order:

In *Section 2*, the needed facts and supplementary results from the book of Mahmudov [15] are given; Hamiltonian function $H$ and argmaximum sets of a set-valued mapping $F$, the LAM, infimal convolution of proper convex functions, conjugate function for Hamiltonian function taken with a minus sign are introduced and the problems for second order DSIs ($P_D$) and DFIs ($P_C$) with initial point constraints are formulated.

In *Section 3* necessary and sufficient conditions of optimality for second order DSIs are formulated, the dual problem for second order DSIs ($P_D$) is constructed. In what follows, we prove that if $\alpha$ and $\alpha^*$ are the values of primary and dual problems, respectively, then $\alpha \geq \alpha^*$ for all feasible solutions. Moreover, if a certain "nondegeneracy condition", that is standard condition of convex analysis on existence of interior point, is satisfied, then the existence of a solution to one of these problems implies the existence of a solution to the other problem, where $\alpha = \alpha^*$, and in the case where $\alpha > -\infty$ the dual problem has a solution. Finally, duality relationship between a pair of optimization problems with initial point constraint established; it is proved that the Euler-Lagrange type adjoint DSI at the same time is a dual relation.

*Section 4* is devoted to duality of discrete-approximate problem ($P_{DA}$). Consequently, by using the first and the second order difference operators and auxiliary set-valued mapping, the problem for second order continuous-time evolution inclusions ($P_C$) is approximated with associated discrete-approximation problem ($P_{DA}$). It is noted that, transition to the problem ($P_{DA}$) requires some special results on calculation of conjugate and support functions, connecting the dual problems ($P_D^*$) and ($P_{DA}^*$) to DSIs ($P_D$) and discrete-approximate ($P_{DA}$) problems, respectively. For construction of the duality for a Mayer problem in its general form, and positive-homogeneity of Hamilton type functions skilfully are used. Here existence of first and second order difference derivatives considerably make difficult to calculate a conjugate Mayer functional and support functions. It is obvious that this method, which is certainly of independent interest from qualitative viewpoint, can play an important role in numerical procedures as well. At the end of this section is considered a Mayer problem with so-called semilinear second order discrete-time inclusions.

In *Section 5* sufficient conditions of optimality and dual problems $(P_C^*)$ for convex DFIs $(P_C)$ are deduced, where establishment of the dual problem is obtained by passing to the formal limit in dual problem for discrete-approximate problem $(P_{DA}^*)$. Besides, it is proven that the Euler-Lagrange type inclusion is a dual relation in continuous problem, too. The considered semilinear problem shows that maximization in the dual problems is realized over the set of solutions of the adjoint equation. In addition, the optimal values in the primary convex $(P_C)$ and the dual concave $(P_C^*)$ problems are equal: $\inf(P_C) = \sup(P_C^*)$. Thus, it is proved that the Euler-Lagrange type adjoint inclusion at the same time is a dual relation, that is a pair of solutions of primary and dual problems satisfies this relation.

In this work we pursue a twofold goal. First, we constructed a dual problem for a discrete-approximate problem to continuous problem $(P_C)$. Second, we use this direct method to establish a dual problem to a continuous Mayer problem. The construction of a dual problem to the latter is implemented by passing to the formal limit as the discrete step tends to zero. Finally, relying to the described method, we believe that within the framework of the idea of this paper a dual problem can be obtained for any higher order differential inclusions.

## 2. Needed Facts and Problem Statement

Further, for the convenience of the reader, all the necessary concepts, definitions of a convex analysis can be found in the book of Mahmudov [15]. Let $\mathbb{R}^n$ be a $n$-dimensional Euclidean space, $\langle x, u \rangle$ be an inner product of elements $x, y \in \mathbb{R}^n$, and $(x, y)$ be a pair of $x, y$. Assume that $G: \mathbb{R}^n \times \mathbb{R}^n \rightrightarrows \mathbb{R}^n$ is a set-valued mapping from $\mathbb{R}^{2n} = \mathbb{R}^n \times \mathbb{R}^n$ into the set of subsets of $\mathbb{R}^n$. Then $G: \mathbb{R}^{2n} \rightrightarrows \mathbb{R}^n$ is convex if its $\operatorname{gph} G = \{(x, y, z) : z \in G(x, y)\}$ is a convex subset of $\mathbb{R}^{3n}$. The set-valued mapping $G$ is convex closed if its graph is a convex closed set in $\mathbb{R}^{3n}$. The domain of $G$ is denoted by $\operatorname{dom} G$ and is defined as follows $\operatorname{dom} G = \{(x, y) : G(x, y) \neq \varnothing\}$. $G$ is convex-valued if $G(x, y)$ is a convex set for each $(x, y) \in \operatorname{dom} G$. Let us introduce the Hamiltonian function and argmaximum set for a set-valued mapping $G$

$$H_G(x, y, z^*) = \sup_z \{\langle z, z^* \rangle : z \in G(x, y)\}, \; z^* \in R^n,$$

$$G_A(x, y, z^*) = \{z \in G(x, y) : \langle z, z^* \rangle = H_G(x, y, z^*)\},$$

respectively. For convex $G$ we set $H_G(x, y, z^*) = -\infty$ if $G(x, y) = \varnothing$.

As usual, $W_A(x^*)$ is a support function of the set $A \subset \mathbb{R}^n$, i.e.,

$$W_A(x^*) = \sup_x \{\langle x, x^* \rangle : x \in A\}, \; x^* \in \mathbb{R}^n.$$

Let $\operatorname{int} A$ be the interior of the set $A \subset \mathbb{R}^{3n}$ and $\operatorname{ri} A$ be the relative interior of the set $A$, i.e. the set of interior points of $A$ with respect to its affine hull $\operatorname{Aff} A$.

The convex cone $K_A(w_0)$, $w_0 = (x_0, y_0, z_0)$ is called the cone of tangent directions at a point $w_0 \in A$ to the set $A$ if from $\bar{w} = (\bar{x}, \bar{y}, \bar{z}) \in K_A(w_0)$ it follows that $\bar{w}$ is a tangent vector to the set $A$ at point $w_0 \in A$, i.e., there exists such function $\eta : \mathbb{R}^1 \to \mathbb{R}^{3n}$ that $w_0 + \lambda \bar{w} + \eta(\lambda) \in A$ for sufficiently small $\lambda > 0$ and $\lambda^{-1}\eta(\lambda) \to 0$, as $\lambda \downarrow 0$.

A function $\varphi$ is called a proper function if it does not assume the value $-\infty$ and is not identically equal to $+\infty$. Obviously, $\varphi$ is proper if and only if $\mathrm{dom}\,\varphi \neq \varnothing$ and $\varphi(x, y)$ is finite for $(x, y) \in \mathrm{dom}\,\varphi = \{(x, y) : \varphi(x, y) < +\infty\}$.

In general, for a set-valued mapping $G$ a set-valued mapping $G^*(\cdot, x, y, z) : \mathbb{R}^n \rightrightarrows \mathbb{R}^{2n}$ defined by
$$G^*\left(z^*; (x, y, z)\right) := \left\{(x^*, y^*) : (x^*, y^*, -z^*) \in K^*_{gphG}(x, y, z)\right\},$$
is called the LAM to a set-valued $G$ at a point $(x, y, z) \in gphG$, where $K^*_{gphG}(x, y, z)$ is the dual to the cone of tangent directions $K_{gphG}(x, y, z)$. We provide another definition of LAM to mapping $G$ which is more relevant for further development
$$G^*(z^*; (x, y, z)) := \left\{(x^*, y^*) : H_G(x_1, y_1, z^*) - H_G(x, y, z^*) \leq \langle x^*, x_1 - x \rangle \right.$$
$$\left. + \langle y^*, y_1 - y \rangle, \forall (x_1, y_1) \in \mathbb{R}^{2n}\right\}, \quad (x, y, z) \in gphG, \quad z \in G_A(x, y, z^*).$$

Clearly, for the convex mapping the Hamiltonian $H(\cdot, \cdot, z^*)$ is concave and the latter and previous definitions of LAMs coincide.

**Definition 2.1** A function $\varphi(x, y)$ is said to be a closure if its epigraph $\mathrm{epi}\,\varphi = \{(x^0, x, y) : x^0 \geq \varphi(x, y)\}$ is a closed set.

**Definition 2.2** The function $\varphi^*(x^*, y^*) = \sup_{x, y} \{\langle x, x^* \rangle + \langle y, y^* \rangle - \varphi(x, y)\}$ is called the conjugate of $\varphi$. It is clear to see that the conjugate function is closed and convex.

Let us denote
$$M_G(x^*, y^*, z^*) = \inf_{x, y, z} \left\{\langle x, x^* \rangle + \langle y, y^* \rangle - \langle z, z^* \rangle : (x, y, z) \in \mathrm{gph}G\right\},$$
that is, for every $(x, y) \in \mathbb{R}^{2n}$
$$M_G(x^*, y^*, z^*) \leq \langle x, x^* \rangle + \langle y, y^* \rangle - H_G(x, y, z^*).$$
It is clear that the function
$$M_G(x^*, y^*, z^*) = \inf_{x, y} \left\{\langle x, x^* \rangle + \langle y, y^* \rangle - H_G(x, y, z^*)\right\}$$
is a support function taken with a minus sign. Besides, it follows that for a fixed $z^*$
$$M_G(x^*, y^*, z^*) = -\left[H_G(\cdot, \cdot, z^*)\right]^*(x^*, y^*)$$
that is, $M_G$ is the conjugate function for $H_G(\cdot, \cdot, z^*)$ taken with a minus sign.

**Definition 2.3** We recall that the operation of infimal convolution $\oplus$ of functions $f_i, i = 1, \ldots, k$ is defined as follows
$$(f_1 \oplus \cdots \oplus f_k)(u) = \inf\left\{f_1(u^1) + \ldots + f_k(u^k) : u^1 + \ldots + u^k = u\right\}, u^i \in \mathbb{R}^m, i = 1, \ldots, k.$$

The infimal convolution $(f_1 \oplus \cdots \oplus f_k)$ is said to be exact provided the infimum above is attained for every $u \in \mathbb{R}^m$. One has $\text{dom}(f_1 \oplus \cdots \oplus f_k) = \sum_{i=1}^{i=k} \text{dom} f_i$. Besides for a proper convex closed functions $f_i, i = 1,...,k$ their infimal convolution $(f_1 \oplus \cdots \oplus f_1)$ is convex and closed (but not necessarily proper). If $f_i, i = 1,...,k$ are functions not identically equal to $+\infty$, then $(f_1 \oplus \cdots \oplus f_k)^* = \sum_{i=1}^{i=k} f_i^*$. Thus, the conjugate of infimal convolution is the sum of the conjugates and this holds without any requirement on the convex functions. The operations + and $\oplus$ are thus dual to each other with respect to taking conjugates.

In Section 5 we deal with the Mayer problem for $(P_C)$ type of the evolution DFIs:

$$\text{infimum } \varphi(x(1), x'(1)), \tag{1}$$

$(P_C)$
$$x''(t) \in F(x(t), x'(t), t), \text{ a.e. } t \in [0,1], \tag{2}$$

$$x(0) \in Q_0, \; x'(0) \in Q_1. \tag{3}$$

Here $F(\cdot, t): \mathbb{R}^n \rightrightarrows \mathbb{R}^n$ is a time dependent set-valued mapping, $\varphi$ is continuous $\varphi: \mathbb{R}^{2n} \to \mathbb{R}^1$, $Q_i \subseteq \mathbb{R}^n (i = 0,1)$ are nonempty subsets. The problem is to find an arc $\tilde{x}(\cdot)$ of the problem (1) – (3) satisfying (2) almost everywhere (a.e.) on $[0,1]$ and the initial-point constraints (3) on $[0,1]$ that minimizes the Mayer functional $\varphi(x(1), x'(1))$. We label this problem as $(P_C)$. Here, a feasible trajectory $x(\cdot)$ is understood to be an absolutely continuous function on a time interval $[0,1]$ together with the first order derivatives for which $x''(\cdot) \in L_1^n([0,1])$. Obviously, such class of functions is a Banach space, endowed with the different equivalent norms.

Section 3 is concerned with the following second order discrete model labelled as $(P_D)$:

$$\text{infimum } \varphi(x_{N-1}, x_N), \tag{4}$$

$(P_D)$
$$x_{t+2} \in F(x_t, x_{t+1}, t), \; t = 0,...,N-2, \tag{5}$$

$$x_0 \in Q_0, \; x_1 \in Q_1. \tag{6}$$

A sequence $\{x_t\}_{t=0}^N = \{x_t : t = 0, 1, ..., N\}$ is called a feasible trajectory for the stated problem (4)-(6). It is required find a solution $\{\tilde{x}_t\}_{t=0}^N$ to a problem $(P_D)$ for the second discrete-time problem, satisfying (5), (6) and minimizing $\varphi(x_{N-1}, x_N)$. In what follows, to this end our further strategy is as follows: first to derive necessary and sufficient conditions of optimality for problems $(P_C)$ and $(P_D)$ and then to derive duality results for them.

**Definition 2.4** Let us say that for the convex problem (4) – (6) the nondegeneracy condition is satisfied if for points $x_t \in \mathbb{R}^n$, one of the following cases is fulfilled:

**(i)** $(x_t, x_{t+1}, x_{t+2}) \in \text{ri gph } F(\cdot, t), x_0 \in \text{ri } Q_0; x_1 \in \text{ri } Q_1, (x_{N-1}, x_N) \in \text{ridom} \varphi, t = 0,...,N-2$,

**(ii)** $(x_t, x_{t+1}, x_{t+2}) \in \text{int gph } F(\cdot, t), \; t = 0,...,N-2, \; x_0 \in \text{int } Q_0; \; x_1 \in \text{int } Q_1$ (with the possible exception of one fixed $t$) and $\varphi$ is continuous at $(x_{N-1}, x_N)$. It follows from the nondegeneracy condition that if $\{\tilde{x}_t\}_{t=0}^N$ is the optimal trajectory in the problem (4)-(6), then the cones of tangent directions $K_{gphG(\cdot,t)}(\tilde{x}_t,\tilde{x}_{t+1},\tilde{x}_{t+2})$ are not separable and consequently, the condition of Theorem 3.2 [15, p.98] is satisfied.

### 3. Infimal Convolution and Duality for a second order DSIs

At first we consider the convex problem (4)-(6). Let us introduce a vector $u = (x_0, x_1, ..., x_N) \in \mathbb{R}^m$, $m = n(N+1)$ and define in the space $\mathbb{R}^{n(N+1)}$ the following convex sets

$$M_t = \{u = (x_0,...,x_N) : (x_t, x_{t+1}, x_{t+2}) \in gph\, F(\cdot,t)\};$$

$$D_0 = \{u = (x_0,...,x_N) : x_0 \in Q_0\},\ D_1 = \{u = (x_0,...,x_N) : x_1 \in Q_1\}.$$

Now, denoting $f(u) = \varphi(x_{N-1}, x_N)$ we will reduce this problem to the problem with geometric constraints. Indeed, it can easily be seen that our basic problem (4)-(6) is equivalent to the following one

$$\text{minimize } f(u) \text{ subject to } A = \left(\bigcap_{t=0}^{N-2} M_t\right) \cap D_0 \cap D_1, \tag{7}$$

where $A$ is a convex set.

In the sense of the terminology of second order DSI [18,19,21,22] we are ready to give the necessary and sufficient conditions for the problem (2.1)-(2.3), which will play an important role in the next investigations.

**Theorem 3.1** Let $F(\cdot,t)$ be an evolution convex set-valued mapping and $Q_i \subseteq \mathbb{R}^n (i=0,1)$ be nonempty convex sets. Besides, let $\{\tilde{x}_t\}_{t=0}^N$ be an optimal trajectory to the second order discrete-time problem with ($P_D$). Then there exist vectors $x_t^*, x_N^*, \mu_t^*$, $t=0,...,N-1$ and a number $\lambda \in \{0,1\}$ simultaneously not all zero, such that the adjoint Euler-Lagrange type inclusions (*i*), and transversality conditions (*ii*), (*iii*) hold:

(i) $\left(x_t^* - \mu_t^*, \mu_{t+1}^*\right) \in F^*\left(x_{t+2}^*; (\tilde{x}_t, \tilde{x}_{t+1}, \tilde{x}_{t+2}), t\right),\ t=0,1,...,N-2$,

(ii) $\mu_0^* - x_0^* \in K_{Q_0}^*(\tilde{x}_0);\ -x_1^* \in K_{Q_1}^*(\tilde{x}_1)$,

(iii) $\left(\mu_{N-1}^* - x_{N-1}^*, -x_N^*\right) \in \lambda \partial_{(x,y)} \varphi(\tilde{x}_{N-1}, \tilde{x}_N)$.

Moreover, if the *nondegeneracy condition* is satisfied, then $\lambda = 1$ and these conditions are sufficient for optimality of the trajectory of $\{\tilde{x}_t\}_{t=0}^N$.

*Proof.* In accordance with the nondegeneracy condition, it follows from Theorems 1.10 and 1.11 [15] that

$$K_A^*(\tilde{u}) = \sum_{t=0}^{N-2} K_{M_t}^*(\tilde{u}) + K_{D_0}^*(\tilde{u}) + K_{D_1}^*(\tilde{u}),\ \tilde{u} = (\tilde{x}_0,...,\tilde{x}_N).$$

Then the rest of the proof is the simple modification of proof of Theorem 5.1[19] and so is omitted. □

For construction of duality we need the following result.

**Proposition 3.1** The conjugate function of $f(u) = \varphi(x_{N-1}, x_N)$ is the following function

$$f^*(x_0^*,...,x_N^*) = \varphi^*(x_{N-1}^*, x_N^*);\ x_i^* = 0,\ i = 0,...,N-2.$$

*Proof.* In fact, the proof of proposition inferred immediately from a definition of conjugate functions:

$$f^*(u^*) = \sup_u \{\langle u, u^* \rangle - f(u)\} = \sup_{x_0,\ldots,x_N} \left\{ \sum_{i=0}^{N} \langle x_i, x_i^* \rangle - \varphi(x_{N-1}, x_N) \right\}$$

$$= \sup_{x_0,\ldots,x_N} \left\{ \sum_{i=0}^{N-2} \langle x_i, x_i^* \rangle + \langle x_{N-1}, x_{N-1}^* \rangle + \langle x_N, x_N^* \rangle - \varphi(x_{N-1}, x_N) \right\}$$

$$= \begin{cases} \varphi^*(x_{N-1}^*, x_N^*), & \text{if } x_i^* = 0, \, i = 0, \ldots, N-2, \\ +\infty, & \text{otherwise.} \end{cases} \qquad \square$$

We call the following problem, labelled as $(P_D^*)$, the dual problem to the problem with second order DSIs $(P_D)$

$$(P_D^*) \quad \sup_{\substack{x_t^*, x_N^*, \mu_t^*, \\ t=0,\ldots,N-1}} \left\{ -\varphi^*(\mu_{N-1}^* - x_{N-1}^*, -x_N^*) + \sum_{t=0}^{N-2} M_{F(\cdot,t)}\left(x_t^* - \mu_t^*, \mu_{t+1}^*, x_{t+2}^*\right) - W_{Q_0}(x_0^* - \mu_0^*) - W_{Q_1}(x_1^*) \right\},$$

where $W_{Q_0}$, $W_{Q_1}$ are support functions of the sets $Q_0, Q_1$, respectively.

**Theorem 3.2** If $\alpha$ and $\alpha^*$ are the optimal values of the optimization problem for second order DSIs $(P_D)$ and its dual problem $(P_D^*)$, respectively, then $\alpha \geq \alpha^*$ for all feasible solutions of primary $(P_D)$ and dual $(P_D^*)$ problems. Besides, if the nondegeneracy condition is satisfied, then, the existence of a solution to one of these problems implies the existence of a solution to another, where $\alpha = \alpha^*$ and in the case $\alpha > -\infty$ the dual problem $(P_D^*)$ has a solution.

*Proof.* It is known from convex analysis that the operations of addition and infimal convolution of convex functions are dual to each other [14]. By this result, if there exists a point $u^0 \in A$, where $f$ is continuous ($f$ is continuous on ri dom$f$, however, $f$ may have a point of discontinuity in its boundary), the optimal value of problem (7) is

$$\inf_{u \in A} f(u) = \inf \{f(u) + \delta_A(u)\} = -\sup \{-f(u) - \delta_A(u)\} = -\sup \{\langle u, 0 \rangle - [f(u) + \delta_A(u)]\}$$

$$= -(f + \delta_A)^*(0) = -(f^* \oplus \delta_A^*)(0) = -\inf \{f^*(u^*) + \delta_A^*(-u^*)\} = \sup \{-f^*(u^*) - \delta_A^*(-u^*)\},$$

where $\delta_A(\cdot)$ is the indicator function of $A$, i.e.,

$$\delta_A(u) = \begin{cases} 0, & u \in A, \\ +\infty, & u \notin A. \end{cases}$$

Note that the nondegeneracy condition guarantees that a point $u^0 \in A$ having this property exists.

In general, it can be noticed that $(f + \delta_A)^*(0) \leq (f^* \oplus \delta_A^*)(0)$ and so

$$\inf_{u \in A} f(u) \geq \sup \{-f^*(u^*) - \delta_A^*(-u^*)\}.$$

Then it is reasonable to announce that the dual problem to the primary problem (7) has the form

$$\sup \{-f^*(u^*) - \delta_A^*(-u^*)\}. \tag{8}$$

In addition, if the value of the problem (7) is finite, then the supremum in the problem (8) is attained for all $u^*$. Recall that, the indicator function of intersection is a sum of indicator functions, that is, $\delta_A = \sum_{t=0}^{N-2} \delta_{M_t} + \delta_{D_0} + \delta_{D_1}$. Then by the reminded above duality theorem, for we have

$$\delta_A^*(-u^*) \leq \inf\left\{\sum_{t=0}^{N-2}\delta_{M_t}^*(-u^*(t)) + \delta_{D_0}^*(-\bar{u}^*(0)) + \delta_{D_1}^*(-\bar{u}^*(1)) : \sum_{i=0}^{N-2}u^*(i) + \bar{u}^*(0) + \bar{u}^*(1) = u^*\right\}, \quad (9)$$

where $u^*(i) = (x_0^*(i),...,x_N^*(i))$, $i = 0,...,N-2$ and $\bar{u}^*(i) = (\bar{x}_0^*(i),...,\bar{x}_N^*(i))$, $i = 0,1$.

In addition, as is well known, the conjugate function of the indicator function of a convex set is the support function of this set, and by Theorem 1.25 [15] the converse assertion is true if the considered set is closed. Then we deduce that

$$\delta_{M_t}^*(-u^*(t)) = \begin{cases} -\inf_{(x_t(t),x_{t+1}(t),x_{t+2}(t))\in \text{gph}F}\left[\langle x_t(t), x_t^*(t)\rangle + \langle x_{t+1}(t), x_{t+1}^*(t)\rangle + \langle x_{t+2}(t), x_{t+2}^*(t)\rangle\right], x_i^*(t) = 0, \\ \hfill i \neq t, t+1, t+2, \\ +\infty, \hfill \text{otherwise,} \end{cases}$$

$$t = 0,...,N-2; \quad (10)$$

$$\delta_{D_0}^*(-\bar{u}^*(0)) = \begin{cases} \sup\left[\langle x_0(0), -\bar{x}_0^*(0)\rangle : x_0(0) \in Q_0\right], & \text{if } \bar{x}_i^*(0) = 0, i \neq 0, \\ +\infty, & \text{otherwise,} \end{cases}$$

$$\delta_{D_1}^*(-\bar{u}^*(1)) = \begin{cases} \sup\left[\langle x_1(1), -\bar{x}_1^*(1)\rangle : x_1(1) \in Q_1\right], & \text{if } \bar{x}_i^*(1) = 0, i \neq 1, \\ +\infty, & \text{otherwise.} \end{cases}$$

(11)

Further, from the relationships (9)-(11) and the formula $f^*(u^*) = \varphi^*(x_{N-1}^*, x_N^*)$ of Proposition 3.1, where $x_i^* = 0, i \neq N-1, N$, with the preceding notations, we conclude that

$$\sup\{-f^*(u^*) - \delta_A^*(-u^*)\} = \sup\left\{-\varphi^*(x_{N-1}^*, x_N^*) + \sum_{t=0}^{N-2} M_{F(\cdot,t)}\left(x_t^*(t), x_{t+1}^*(t), -x_{t+2}^*(t)\right)\right.$$

$$-W_{Q_0}(-\bar{x}_0^*(0)) - W_{Q_1}(-\bar{x}_1^*(1)) : x_0^*(0) + \bar{x}_0^*(0) = 0, \ x_1^*(0) + x_1^*(1) + \bar{x}_1^*(1) = 0,$$

$$x_t^*(t) + x_t^*(t-1) + x_t^*(t-2) = 0 \ (t = 2,...,N-2),$$

$$\left. x_{N-1}^*(N-2) + x_{N-1}^*(N-3) = x_{N-1}^*, \ x_N^*(N-2) = x_N^*\right\}, \quad (12)$$

where the supremum is attained, if $\alpha > -\infty$. For further convenience, we denote $x_{t+1}^*(t) \equiv \mu_{t+1}^*$, $-x_{t+2}^*(t) \equiv x_{t+2}^*$, $t = 0,...,N-2$. Then if we take $\bar{x}_0^*(0) = \mu_0^* - x_0^*$, $\bar{x}_1^*(1) = -x_1^*$, we can generalize the relation $x_t^*(t) + x_t^*(t-1) + x_t^*(t-2) = 0$ in (12) to the case $t = 0,1$. Thus, taking into account these new designations, the right hand side of the relation (12) has the form $(P_D^*)$.

□

In fact, in the next theorem is proved that the Euler-Lagrange type inclusions (i) of Theorem 3.1 is a dual relation for the pairs of primal ($P_D$) and dual ($P_D^*$) problems.

**Theorem 3.3** Suppose that $\{\tilde{x}_t\}_{t=0}^N$ is an optimal solution to primary problem ($P_D$) and that the nondegeneracy condition is satisfied. Besides, suppose that $\varphi$ is proper convex closed function, $Q_0, Q_1$ and $F(\cdot,t)$ are convex closed sets and set-valued mapping, respectively. Then the family of vectors $\{x_t^*, \mu_t^*\}_{t=0}^N$ is an optimal solution to the dual problem $(P_D^*)$ if and only if the adjoint Euler-Lagrange type inclusions (*i*), and transversality conditions (*ii*), (*iii*) of Theorem 3.1 are satisfied.

*Proof.* Let $\{\tilde{x}_t\}_{t=0}^N$ and $\{x_t^*, \mu_t^*\}_{t=0}^N$ be the optimal solutions of problems $(P_D)$ and $(P_D^*)$, respectively. Prove that the conditions (*i*)-(*iii*) of Theorem 3.1 is satisfied. Since the problems $(P_D)$ and (7) are equivalent, then we have $\tilde{u} = (\tilde{x}_0, ..., \tilde{x}_N) \in \{u : f(u) + \delta_A(u) = \alpha\}$. It follows that, $0 \in \partial(f(\tilde{u}) + \delta_A(\tilde{u}))$ or $\tilde{u} \in \partial(f + \delta_A)^*(0) = \partial(f^* \oplus \delta_A^*)(0)$. Then in view of Theorem 6.6.5 [14] on the calculation of the subdifferential of the infimal convolution, we derive

$$\partial(f^* \oplus \delta_A^*)(0) = \partial f^*(u^*) \bigcap \partial \delta_A^*(-u^*).$$

Then the latter inclusion give us

$$\tilde{u} \in \partial f^*(u^*) \bigcap \partial \delta_A^*(-u^*) \neq \varnothing. \tag{13}$$

Therefore, by Proposition 6.7.8 [14], it can be conclude that $u^*$ is a solution to the maximization problem $\sup\{-f^*(u^*) - \delta_A^*(-u^*)\}$. As a consequence, from the nondegeneracy condition and Theorem 6.6.6 [14], we have

$$\partial \delta_A^*(-u^*) = \partial\left[\bigoplus_{t=0}^{N-2} \delta_{M_t}^* \oplus \delta_{D_0}^* \oplus \delta_{D_1}^*\right](-u^*), \quad \sum_{t=0}^{N-2} u^*(t) + \bar{u}^*(0) + \bar{u}^*(1) = u^*,$$

and consequently

$$\partial \delta_A^*(-u^*) = \bigcap_{t=0}^{N-2} \partial \delta_{M_t}^*(u^*(t)) \bigcap \partial \delta_{D_0}^*(\bar{u}^*(0)) \bigcap \partial \delta_{D_1}^*(\bar{u}^*(1)). \tag{14}$$

Then taking into account (13) and (14) we have

$$\tilde{u} \in \partial \delta_{M_t}^*(u^*(t))\, (t=0,...,N-2),\ \tilde{u} \in \partial \delta_{D_0}^*(\bar{u}^*(0)),\ \tilde{u} \in \partial \delta_{D_1}^*(\bar{u}^*(1)),\ \tilde{u} \in \partial f^*(u^*).$$

Now it can be easily seen that since $M_t(t=0,...,N-2), D_0, D_1$ are closed sets, the indicator functions $\delta_{M_t}, \delta_{D_0}, \delta_{D_1}$ are proper convex closed functions. Therefore, the last relationship is equivalent to

$$u^*(t) \in \partial \delta_{M_t}(\tilde{u}),\ t=0,...,N-2;\ \bar{u}^*(0) \in \partial \delta_{D_0}(\tilde{u}),\ \bar{u}^*(1) \in \partial \delta_{D_1}(\tilde{u}),\ u^* \in \partial f(\tilde{u}). \tag{15}$$

Further, it should be noted that on the one hand

$$\partial \delta_{M_t}(\tilde{u}) = -K_{M_t}^*(\tilde{u}),\ t=0,...,N-2;\ \partial \delta_{D_0}(\tilde{u}) = -K_{D_0}^*(\tilde{u}),\ \partial \delta_{D_1}(\tilde{u}) = -K_{D_1}^*(\tilde{u})$$

and on the other hand

$$K_{M_t}^*(\tilde{u}) = \left\{u^*(t) : \left(x_t^*(t), x_{t+1}^*(t), x_{t+2}^*(t)\right) \in K_{gphF(\cdot,t)}^*\left(\tilde{x}_t, \tilde{x}_{t+1}, \tilde{x}_{t+2}\right):\ x_i^*(t) = 0, i \neq t, t+1, t+2\right\},$$

$$K_{D_0}^*(\tilde{u}) = \left\{u^*(t) : \bar{x}_0^*(0) \in K_{D_0}^*\left(\bar{x}_i(0)\right):\ \bar{x}_i^*(0) = 0, i \neq 0\right\},$$

$$K_{D_1}^*(\tilde{u}) = \left\{u^*(t) : \bar{x}_1^*(1) \in K_{D_1}^*\left(\bar{x}_i(1)\right):\ x_i^*(1) = 0, i \neq 1\right\}.$$

Therefore, inclusions (15) give

$$\left(x_t^*(t), x_{t+1}^*(t), x_{t+2}^*(t)\right) \in K_{gphF(\cdot,t)}^*\left(\tilde{x}_t, \tilde{x}_{t+1}, \tilde{x}_{t+2}\right),\ i=0,...,N-2,$$

$$\bar{x}_0^*(0) \in K_{Q_0}^*\left(\bar{x}_0(0)\right),\ \bar{x}_1^*(1) \in K_{Q_1}^*\left(\bar{x}_1(1)\right),\ \left(x_{N-1}^*, x_N^*\right) \in \partial_{(x,y)}\varphi(\tilde{x}_{N-1}, \tilde{x}_N), \tag{16}$$

where the vectors $u^*(t)\,(t=0,...,N-2), \bar{u}^*(0), \bar{u}^*(1)$ satisfy the condition $\sum_{t=0}^{N-2} u^*(t) + \bar{u}^*(0) + \bar{u}^*(1) = u^*$. Then, taking into account the above accepted designations $x_{t+1}^*(t) \equiv \mu_{t+1}^*,\ -x_{t+2}^*(t)$

$\equiv x_{t+2}^{*}$, $t = 0,...,N-2$, $\bar{x}_{0}^{*}(0) = \mu_{0}^{*} - x_{0}^{*}$, $\bar{x}_{1}^{*}(1) = -x_{1}^{*}$, where the collection of vectors $x_{t}^{*}, \mu_{t}^{*}$, $t = 0,...,N$ is a solution to the dual problem, from (16), by applying the definition of LAM we can derive that the conditions (*i*)-(*iii*) of Theorem 3.1 are satisfied.

Now we shall prove the converse assertion. Let for a family of vectors $\{x_{t}^{*}, \mu_{t}^{*}\}_{t=0}^{N}$ the conditions (*i*)-(*iii*) of Theorem 3.1 be satisfied, where $\{\tilde{x}_{t}\}_{t=0}^{N}$ be an optimal solution to the primary problem ($P_D$). We should prove that $\{x_{t}^{*}, \mu_{t}^{*}\}_{t=0}^{N}$ is a solution to the dual problem ($P_{D}^{*}$). It is easy to see that by Lemma 2.6 [15] the adjoint Euler-Lagrange type inclusions

$$(x_{t}^{*} - \mu_{t}^{*}, \mu_{t+1}^{*}) \in F^{*}(x_{t+2}^{*}; (\tilde{x}_{t}, \tilde{x}_{t+1}, \tilde{x}_{t+2}), t), \quad t = 0,1,...,N-2$$

and conditions

$$M_{F(\cdot,t)}(x_{t}^{*} - \mu_{t}^{*}, \mu_{t+1}^{*}, x_{t+2}^{*}) = \langle \tilde{x}_{t}, x_{t}^{*} - \mu_{t}^{*} \rangle + \langle \tilde{x}_{t+1}, \mu_{t+1}^{*} \rangle - H_{F}(\tilde{x}_{t}, \tilde{x}_{t+1}, x_{t+2}^{*}), \quad t = 0,...,N-2$$

are equivalent. Then in view of the condition $\tilde{x}_{t+2} \in F_{A}(\tilde{x}_{t}, \tilde{x}_{t+1}, x_{t+2}^{*})$ or, another speaking, in view of $\langle \tilde{x}_{t+2}, x_{t+2}^{*} \rangle = H_{F}(\tilde{x}_{t}, \tilde{x}_{t+1}, x_{t+2}^{*})$, from Theorem 3.1 we get

$$M_{F(\cdot,t)}(x_{t}^{*} - \mu_{t}^{*}, \mu_{t+1}^{*}, x_{t+2}^{*}) = \langle \tilde{x}_{t}, x_{t}^{*} - \mu_{t}^{*} \rangle + \langle \tilde{x}_{t+1}, \mu_{t+1}^{*} \rangle - \langle \tilde{x}_{t+2}, x_{t+2}^{*} \rangle, \quad t = 0,...,N-2. \quad (17)$$

Since nondegeneracy condition is satisfied and $\lambda = 1$, by Theorem 1.27 [15], the inclusion $(\mu_{N-1}^{*} - x_{N-1}^{*}, -x_{N}^{*}) \in \lambda \partial_{(x,y)} \varphi(\tilde{x}_{N-1}, \tilde{x}_{N})$ is equivalent to

$$\varphi^{*}(\mu_{N-1}^{*} - x_{N-1}^{*}, -x_{N}^{*}) = \langle \tilde{x}_{N-1}, \mu_{N-1}^{*} - x_{N-1}^{*} \rangle - \langle \tilde{x}_{N}, x_{N}^{*} \rangle - \varphi(\tilde{x}_{N-1}, \tilde{x}_{N}). \quad (18)$$

On the other hand, it is not hard to see that the inclusions $\mu_{0}^{*} - x_{0}^{*} \in K_{Q_{0}}^{*}(\tilde{x}_{0})$; $-x_{1}^{*} \in K_{Q_{1}}^{*}(\tilde{x}_{1})$ imply

$$W_{Q_{0}}(x_{0}^{*} - \mu_{0}^{*}) = -\langle \mu_{0}^{*} - x_{0}^{*}, \tilde{x}_{0} \rangle; \quad W_{Q_{1}}(x_{1}^{*}) = \langle x_{1}^{*}, \tilde{x}_{1} \rangle. \quad (19)$$

As a result of summation of the equalities (17)-(19) we have

$$-\varphi^{*}(\mu_{N-1}^{*} - x_{N-1}^{*}, -x_{N-1}^{*}, -x_{N}^{*}) + \sum_{t=0}^{N-2} M_{F(\cdot,t)}(x_{t}^{*} - \mu_{t}^{*}, \mu_{t+1}^{*}, x_{t+2}^{*}) - W_{Q_{0}}(x_{0}^{*} - \mu_{0}^{*})$$
$$-W_{Q_{1}}(x_{1}^{*}) = \varphi(\tilde{x}_{N-1}, \tilde{x}_{N}) + \langle \tilde{x}_{N}, x_{N}^{*} \rangle - \langle \tilde{x}_{N-1}, \mu_{N-1}^{*} - x_{N-1}^{*} \rangle + \langle \mu_{0}^{*} - x_{0}^{*}, \tilde{x}_{0} \rangle$$
$$- \langle x_{1}^{*}, \tilde{x}_{1} \rangle + \sum_{t=0}^{N-2} \left[ \langle \tilde{x}_{t}, x_{t}^{*} - \mu_{t}^{*} \rangle + \langle \tilde{x}_{t+1}, \mu_{t+1}^{*} \rangle - \langle \tilde{x}_{t+2}, x_{t+2}^{*} \rangle \right]. \quad (20)$$

Now we can simplify the last term in the right hand side of (20) as follows

$$\sum_{t=0}^{N-2} \left[ \langle \tilde{x}_{t}, x_{t}^{*} - \mu_{t}^{*} \rangle + \langle \tilde{x}_{t+1}, \mu_{t+1}^{*} \rangle - \langle \tilde{x}_{t+2}, x_{t+2}^{*} \rangle \right] = \langle \tilde{x}_{0}, x_{0}^{*} \rangle + \langle \tilde{x}_{1}, x_{1}^{*} \rangle + \sum_{t=2}^{N-2} \langle \tilde{x}_{t}, x_{t}^{*} \rangle$$
$$- \langle \tilde{x}_{0}, \mu_{0}^{*} \rangle - \sum_{t=1}^{N-2} \langle \tilde{x}_{t}, \mu_{t}^{*} \rangle + \sum_{t=1}^{N-1} \langle \tilde{x}_{t}, \mu_{t}^{*} \rangle - \sum_{t=2}^{N} \langle \tilde{x}_{t}, x_{t}^{*} \rangle = \langle \tilde{x}_{0}, x_{0}^{*} \rangle + \langle \tilde{x}_{1}, x_{1}^{*} \rangle \quad (21)$$
$$- \langle \tilde{x}_{N-1}, x_{N-1}^{*} \rangle - \langle \tilde{x}_{N}, x_{N}^{*} \rangle - \langle \tilde{x}_{0}, \mu_{0}^{*} \rangle + \langle \tilde{x}_{N-1}, \mu_{N-1}^{*} \rangle.$$

Finally, introducing (21) in (20) we have

$$-\varphi^{*}(\mu_{N-1}^{*} - x_{N-1}^{*}, -x_{N-1}^{*}, -x_{N}^{*}) + \sum_{t=0}^{N-2} M_{F(\cdot,t)}(x_{t}^{*} - \mu_{t}^{*}, \mu_{t+1}^{*}, x_{t+2}^{*}) - W_{Q_{0}}(x_{0}^{*} - \mu_{0}^{*})$$
$$-W_{Q_{1}}(x_{1}^{*}) = \varphi(\tilde{x}_{N-1}, \tilde{x}_{N}) + \langle \tilde{x}_{N}, x_{N}^{*} \rangle - \langle \tilde{x}_{N-1}, \mu_{N-1}^{*} - x_{N-1}^{*} \rangle + \langle \mu_{0}^{*} - x_{0}^{*}, \tilde{x}_{0} \rangle - \langle \tilde{x}_{1}, x_{1}^{*} \rangle$$

$$+\langle \tilde{x}_0, x_0^*\rangle + \langle \tilde{x}_1, x_1^*\rangle - \langle \tilde{x}_{N-1}, x_{N-1}^*\rangle - \langle \tilde{x}_N, x_N^*\rangle - \langle \tilde{x}_0, \mu_0^*\rangle + \langle \tilde{x}_{N-1}, \mu_{N-1}^*\rangle = \varphi(\tilde{x}_{N-1}, \tilde{x}_N).$$

Thus, taking the supremum in the left hand side of the latter relation on the set of arbitrary vectors $x_t^*, \mu_t^*$, $t=0,...,N$, we may conclude that $\alpha^* \geq \alpha$. Comparing this with the opposite inequality $\alpha \geq \alpha^*$ (Theorem 3.2), we get $\alpha = \alpha^*$. Consequently, the collection of vectors $x_t^*, \mu_t^*, t=0,...,N$ satisfying the conditions (*i*)-(*iii*) of Theorem 3.1 is an optimal solution to the dual problem $(P_D^*)$. The proof is completed. □

**Example 3.1** Consider the following optimal control problem with DSIs :

$$\text{infimum } \varphi(x_{N-1}, x_N),$$
$$\text{subject to } \quad x_{t+2} = A_0 x_t + A_1 x_{t+1} + Bu_t, \; u_t \in U, \quad (22)$$
$$x_0 \in Q_0, \; x_1 \in Q_1, \; t=0,...,N-2,$$

where $A_0, A_1$ and $B$ are $n \times n$ and $n \times r$ matrices, respectively, $\varphi$ proper convex function, $U \subseteq \mathbb{R}^r$ nonempty convex set, $Q_i \subset \mathbb{R}^n (i=0,1)$ convex sets. The problem is to find a sequence $\{\tilde{u}_t\}_{t=0}^{N-2}$ of controlling parameters $\tilde{u}_t \in U$ such that the corresponding trajectory $\{\tilde{x}_t\}_{t=0}^{N}$ in (22) minimizes $\varphi(x_{N-1}, x_N)$. Before all we introduce a set-valued mapping of the form $F(x_t, x_{t+1}) = A_0 x_t + A_1 x_{t+1} + BU$, $t=0,...,N-2$. Then it is easy to see that

$$M_F(x^*, y^*, z^*) = \inf_{(x,y,z) \in \text{gph}F} \{\langle x, x^*\rangle + \langle y, y^*\rangle - \langle z, z^*\rangle\}$$
$$= \inf_{x,y} \left[\langle x, x^* - A_0^* z^*\rangle + \langle y, y^* - A_1^* z^*\rangle\right] - \sup_{u \in U} \langle u, B^* z^*\rangle$$
$$= \begin{cases} -W_U(B^* z^*), & \text{if } x^* = A_0^* z^*, \; y^* = A_1^* z^*, \\ -\infty, & \text{otherwise.} \end{cases}$$

Therefore,

$$M_F\left(x_t^* - \mu_t^*, \mu_{t+1}^*, x_{t+2}^*\right) = \begin{cases} -W_U(B^* x_{t+2}^*), & \text{if } x_t^* - \mu_t^* = A_0^* x_{t+2}^*, \; \mu_{t+1}^* = A_1^* x_{t+2}^*, \\ -\infty, & \text{otherwise} \end{cases}$$

and

$$M_F\left(x_t^* - \mu_t^*, \mu_{t+1}^*, x_{t+2}^*\right) = \begin{cases} -W_U(B^* x_{t+2}^*), & \text{if } x_t^* = A_0^* x_{t+2}^* + A_1^* x_{t+1}^*, \\ -\infty, & \text{otherwise.} \end{cases}$$

Then according to the problem $(P_D^*)$ the dual problem to problem (22) is

$$\sup_{x_t^*, t=0,...,N} \left\{-\varphi^*(A_1^* x_N^* - x_{N-1}^*, -x_N^*) - \sum_{t=0}^{N-2} W_U(B^* x_{t+2}^*) - W_{Q_0}(x_0^* - A_1^* x_1^*) - W_{Q_1}(x_1^*)\right\},$$
$$x_t^* = A_0^* x_{t+2}^* + A_1^* x_{t+1}^*.$$

Hence, it is interesting to note that supremum here is taken over the set of solutions of the discrete Euler-Lagrange inclusion/equation.

## 4. Construction of Dual Problem for a Discrete-approximate Problem

In this section first of all we should construct the dual problem to discrete-approximate problem, associated to continuous problem $(P_C)$; assume that $\delta$ is a step on the $t$-axis and $x(t) \equiv x_\delta(t)$ is a grid function on a uniform grid on $[0,1]$. We introduce the following first and second order difference operators (forward and backward difference approximations) $\Delta_+ x(t) \equiv \Delta x(t)$ and $\Delta_- x(t)$, $t = 0, \delta, \ldots, 1 - 2\delta$

$$\Delta x(t) = \frac{1}{\delta}[x(t+\delta) - x(t)], \quad \Delta_- x(t) = \frac{1}{\delta}[x(t) - x(t-\delta)], \quad \Delta^2 x(t) = \frac{1}{\delta}[\Delta x(t+\delta) - \Delta x(t)]$$

and associate with the problem ($P_C$) the following second order discrete-approximate evolution problem

$$\text{minimize } \varphi(x(1-\delta), \Delta_- x(1)),$$
$$\Delta^2 x(t) \in F(x(t), \Delta x(t), t), \ t = 0, \delta, \ldots, 1 - 2\delta, \tag{23}$$
$$x(0) \in Q_0, \ \Delta x(0) \in Q_1.$$

We should reduce the problem (23) to a problem of the form (4)-(6) or ($P_D$) with initial point constraints. Introducing a new set-valued mapping $G(\cdot, t) : \mathbb{R}^n \times \mathbb{R}^n \rightrightarrows \mathbb{R}^n$ and functions

$$G(x, y, t) = 2y - x + \delta^2 F(x, (y-x)/\delta, t),$$
$$\Phi(x(1-\delta), x(1)) \equiv \varphi(x(1), \Delta_- x(1)), \quad \hat{Q}_1 = x(0) + \delta Q_1$$

we rewrite the problem (23) in terms of function $\Phi$ and set-valued mapping $G(\cdot, t)$ as follows:

$$\text{minimize } \Phi(x(1-\delta), x(1)),$$
($P_{DA}$) $\quad x(t + 2\delta) \in G(x(t), x(t+\delta), t),$ $\tag{24}$
$$x(0) \in Q_0, \ x(\delta) \in \hat{Q}_1, \ t = 0, \delta, \ldots, 1 - 2\delta.$$

Now, our main problem is to construct a dual problem for the problem (24). According to the dual problem ($P_D^*$) for problem ($P_{DA}$) we have the following dual problem, labelled by ($P_{DA}^*$)

$$(P_{DA}^*) \quad \sup\left\{ -\Phi^*\left(\bar{\mu}^*(1-\delta) - \bar{x}^*(1-\delta), -\bar{x}^*(1)\right) + \sum_{t=0}^{1-2\delta} M_{G(\cdot, t)}\left(\bar{x}^*(t) - \bar{\mu}^*(t), \bar{\mu}^*(t+\delta), \bar{x}^*(t+2\delta)\right) \right.$$
$$\left. - W_{Q_0}(\bar{x}^*(0) - \bar{\mu}^*(0)) - W_{\hat{Q}_1}(\bar{x}^*(\delta)) \right\}.$$

Here we should express the functions $M_{G(\cdot, t)}$ and $\Phi^*$ in terms $M_{F(\cdot, t)}$ and $\varphi^*$, respectively, which plays a central role in our developments in the next stages. We prove the following simple and important additional results.

**Proposition 4.1** Let $F(\cdot, t) : \mathbb{R}^{2n} \rightrightarrows \mathbb{R}^n$ be a convex set-valued mapping and $G(x, y, t) = 2y - x + \delta^2 F(x, (y-x)/\delta, t)$. Then one has

$$M_{G(\cdot, t)}(x^*, y^*, z^*) = \delta^2 M_{F(\cdot, t)}\left(\frac{x^* + y^* - z^*}{\delta^2}, \frac{y^* - 2z^*}{\delta}, z^*\right).$$

*Proof.* Indeed, by definition of a set-valued mapping $G(\cdot, t)$ and function $M_{G(\cdot, t)}$ we can write

$$M_{G(\cdot, t)}(x^*, y^*, z^*) = \inf_{x, y, z}\left\{\langle x, x^*\rangle + \langle y, y^*\rangle - \langle z, z^*\rangle : (x, y, z) \in \text{gph}G(\cdot, t)\right\}$$
$$= \inf\left\{\langle x, x^*\rangle + \langle y, y^*\rangle - \langle z, z^*\rangle : \left(x, \frac{y-x}{\delta}, \frac{z - 2y + x}{\delta^2}\right) \in \text{gph}F(\cdot, t)\right\}.$$

Now, the expression $\langle x, x^*\rangle + \langle y, y^*\rangle - \langle z, z^*\rangle$ in the last braces we should rewrite in a more

relevant form as

$$\langle x, x_1^* \rangle + \left\langle \frac{y-x}{\delta}, y_1^* \right\rangle - \left\langle \frac{z-2y+x}{\delta^2}, z_1^* \right\rangle, \tag{25}$$

where $x_1^*$, $y_1^*$, $z_1^*$ must be determined. To this end, we rewrite this sum in the form

$$\left\langle x, x_1^* - \frac{y_1^*}{\delta} - \frac{z_1^*}{\delta^2} \right\rangle + \left\langle y, \frac{y_1^*}{\delta} + \frac{2z_1^*}{\delta^2} \right\rangle - \left\langle z, \frac{z_1^*}{\delta^2} \right\rangle$$

and compare it with the expression $\langle x, x^* \rangle + \langle y, y^* \rangle - \langle z, z^* \rangle$. As a result, we have

$$x^* = x_1^* - \frac{y_1^*}{\delta} - \frac{z_1^*}{\delta^2}, \quad y^* = \frac{y_1^*}{\delta} + \frac{2z_1^*}{\delta^2}, \quad z^* = \frac{z_1^*}{\delta^2},$$

whence $x_1^* = x^* + y^* - z^*$, $y_1^* = \delta(y^* - 2z^*)$, $z_1^* = \delta^2 z^*$. Substituting these values into (25) we conclude that

$$\langle x, x^* + y^* - z^* \rangle + \left\langle \frac{y-x}{\delta}, \delta(y^* - 2z^*) \right\rangle - \left\langle \frac{z-2y+x}{\delta^2}, \delta^2 z^* \right\rangle.$$

Finally, we have the desired result:

$$M_{G(\cdot,t)}(x^*, y^*, z^*) = \delta^2 \inf \left\{ \left\langle x, \frac{x^* + y^* - z^*}{\delta^2} \right\rangle + \left\langle \frac{y-x}{\delta}, \frac{y^* - 2z^*}{\delta} \right\rangle - \left\langle \frac{z-2y+x}{\delta^2}, z^* \right\rangle \right.$$

$$\left. : \left( x, \frac{y-x}{\delta}, \frac{z-2y+x}{\delta^2} \right) \in \text{gph} F(\cdot, t) \right\} = \delta^2 M_{F(\cdot,t)} \left( \frac{x^* + y^* - z^*}{\delta^2}, \frac{y^* - 2z^*}{\delta}, z^* \right). \quad \square$$

**Proposition 4.2** If $\delta(\delta > 0)$ is a constant real number, then for a function defined by $\Phi(x, y) \equiv \varphi\left( x, \frac{1}{\delta}(y-x) \right)$ the conjugate function is computed as follows:

(1) $\Phi^*(x^*, y^*) = \varphi^*(x^* + y^*, \delta y^*)$;

In addition, for a proper convex function (not necessarily convex) $\varphi$ the following statement is true:

(2) $(x^*, y^*) \in \partial_{(x,y)} \Phi(x, y)$, $(x, y) \in \text{dom}\,\Phi$, (3) $(x^* + y^*, \delta y^*) \in \partial \varphi\left( x, \frac{1}{\delta}(y-x) \right)$.

*Proof.* By definition of conjugate function, we have

$$\Phi^*(x^*, y^*) = \sup_{x,y} \{ \langle x, x^* \rangle + \langle y, y^* \rangle - \Phi(x, y) \}$$

$$= \sup_{x,y} \left\{ \langle x, x^* \rangle + \left\langle \frac{1}{\delta}(y-x), \delta y^* \right\rangle + \langle x, y^* \rangle - \Phi(x, y) \right\}$$

$$= \sup_{x,y} \left\{ \langle x, x^* + y^* \rangle + \left\langle \frac{1}{\delta}(y-x), \delta y^* \right\rangle - \varphi\left( x, \frac{1}{\delta}(y-x) \right) \right\} = \varphi^*(x^* + y^*, \delta y^*). \tag{26}$$

On the other hand, by Theorem 1.27 [15] it is clear to see that $(x^*, y^*) \in \partial_{(x,y)} \Phi(x, y)$ if and only if $\langle x, x^* \rangle + \langle y, y^* \rangle - \Phi(x, y) = \Phi^*(x^*, y^*)$. Then in view of the equality (26) we deduce that

$$\Phi^*(x^*, y^*) = \varphi^*(x^* + y^*, \delta y^*) = \langle x, x^* + y^* \rangle + \left\langle \frac{1}{\delta}(y-x), \delta y^* \right\rangle - \varphi\left( x, \frac{1}{\delta}(y-x) \right)$$

if and only if $(x^* + y^*, \delta y^*) \in \partial\varphi\left(x, \frac{1}{\delta}(y - x)\right)$. The proof is completed. □

Notice that an important generalization of Proposition 4.2 is the following result.

**Lemma 4.1** Let $E$ be an identity matrix of size $n \times n$ and O be $n \times n$ zero matrix. Besides, let $A = \begin{bmatrix} E & O \\ -E/\delta & E/\delta \end{bmatrix}$ be $2n \times 2n$ block matrix and $\Phi(w) = \varphi(Aw)$, $w = (x, y)$. Then the statement (1) of Proposition 4.2 holds.

*Proof.* Indeed, denoting $u = Aw$ and recalling that $A$ is invertible, on the definition of conjugate functions we have

$$\Phi^*(w^*) = \sup_w\{\langle w, w^*\rangle - \varphi(Aw)\} = \sup_u\{\langle A^{-1}u, w^*\rangle - \varphi(u)\} = \sup_u\{\langle u, A^{*-1}w^*\rangle - \varphi(u)\}$$

$= \varphi^*(A^{*-1}w^*)$. Then, it can be easily verified that $A^{*-1} = \begin{bmatrix} E & E \\ O & \delta E \end{bmatrix}$ and so $A^{*-1}w^*$

$= (x^* + y^*, \delta y^*)$. Consequently, $\Phi^*(w^*) = \varphi^*(A^{*-1}w^*) = \varphi^*(x^* + y^*, \delta y^*)$. □

Now we need the useful result found in Lemma 4.2.

**Remark 4.1.** We remark that if we consider a Mayer problem with third order differential inclusions then we should compute the conjugate function of function

$$\Phi(x, y_1, y_2) \equiv \varphi\left(x, \frac{y_1 - x}{\delta}, \frac{y_2 - 2y_1 + x}{\delta^2}\right).$$

In this case it is not hard to check that $\Phi^*(x^*, y_1^*, y_2^*) = \varphi^*(x^* + y_1^* + y_2^*, \delta y_1^* + 2\delta y_2^*, \delta^2 y_2^*)$.

But increasing the "order" of the problem complicates the calculation. Therefore, noticing that the coefficients of the vectors $x^*, \{y_1^*, \delta y_1^*\}, \{y_2^*, 2\delta y_2^*, \delta^2 y_2^*\}$ form a Pascal triangle with binomial coefficients, we can successfully use it for further "higher order" calculations. In just the same way can be calculated the function $M_{G(\cdot,t)}$. Consequently, all results considered in this work can be generalized to the case of higher order problems.

**Lemma 4.2** Let $F(\cdot, t): \mathbb{R}^{2n} \rightrightarrows \mathbb{R}^n$ be a convex set-valued mapping and $G(x, y, t) = 2y - x + \delta^2 F(x, (y-x)/\delta, t)$. Then between the functions $M_{G(\cdot,t)}$ and $M_{F(\cdot,t)}$ there is the following connection

$$M_{G(\cdot,t)}\left(\bar{x}^*(t) - \bar{\mu}^*(t), \bar{\mu}^*(t+\delta), \bar{x}^*(t+2\delta)\right)$$
$$= \delta M_{F(\cdot,t)}\left(\Delta^2 x^*(t) + \Delta_- v^*(t+\delta), v^*(t+\delta), x^*(t+2\delta)\right),$$

where $v^*(t+\delta) = [\mu^*(t+\delta) - 2x^*(t+2\delta)]/\delta$, $\mu^*(t) = \delta\bar{\mu}^*(t)$, $x^*(t) = \delta\bar{x}^*(t)$.

*Proof.* By applying the Proposition 4.1(the function $M_{F(\cdot,t)}$ is positive homogeneous) we have

$$M_{G(\cdot,t)}\left(\bar{x}^*(t) - \bar{\mu}^*(t), \bar{\mu}^*(t+\delta), \bar{x}^*(t+2\delta)\right)$$

$$= \delta M_{F(\cdot,t)}\left(\frac{x^*(t) - \mu^*(t) + \mu^*(t+\delta) - x^*(t+2\delta)}{\delta^2}, \frac{\mu^*(t+\delta) - 2x^*(t+2\delta)}{\delta}, x^*(t+2\delta)\right). \quad (27)$$

Then taking into account the substitution $\mu^*(t+\delta) = \delta v^*(t+\delta) + 2x^*(t+2\delta)$ it turns out that the following important representation is valid

$$\frac{1}{\delta^2}\left[x^*(t) - \mu^*(t) + \mu^*(t+\delta) - x^*(t+2\delta)\right] = \Delta^2 x^*(t) + \Delta_- v^*(t+\delta).$$

Thus, introducing this formula in (27) immediately we have the desired result. □

**Lemma 4.3** Suppose that a function $\Phi(x(1-\delta), x(1)) \equiv \varphi(x(1), \Delta_- x(1))$ is a proper convex function. Then between the conjugate functions $\Phi^*$ and $\varphi^*$ there is the following connection

$$\Phi^*\left(\bar{\mu}^*(1-\delta) - \bar{x}^*(1-\delta), -\bar{x}^*(1)\right) = \varphi^*\left(v^*(1-\delta) + \Delta_- x^*(1), -x^*(1)\right),$$

where the functions $v^*(t)$, $\mu^*(t)$, $x^*(t)$ are defined as in Lemma 4.2.

*Proof.* According to Proposition 4.2

$$\Phi^*\left(\bar{\mu}^*(1-\delta) - \bar{x}^*(1-\delta), -\bar{x}^*(1)\right) = \varphi^*\left(\bar{\mu}^*(1-\delta) - \bar{x}^*(1-\delta) - \bar{x}^*(1), -\delta \bar{x}^*(1)\right)$$

$$= \varphi^*\left(\frac{\delta \bar{\mu}^*(1-\delta) - \delta \bar{x}^*(1-\delta) - \delta \bar{x}^*(1)}{\delta}, -\delta \bar{x}^*(1)\right) = \varphi^*\left(\frac{\mu^*(1-\delta) - x^*(1-\delta) - x^*(1)}{\delta}, -x^*(1)\right)$$

$$= \varphi^*\left(\frac{\mu^*(1-\delta) - 2x^*(1) + x^*(1) - x^*(1-\delta)}{\delta}, -x^*(1)\right) = \varphi^*\left(v^*(1-\delta) + \Delta_- x^*(1), -x^*(1)\right). \quad \square$$

**Proposition 4.3** The support functions $W_{Q_0}, W_{Q_1}, W_{\hat{Q}_1}$ of the sets $Q_0$, $Q_1$, $\hat{Q}_1 = x(0) + \delta Q_1$, respectively, are connected with the following inequality relation

$$W_{Q_0}\left(\bar{x}^*(0) - \bar{\mu}^*(0)\right) + W_{\hat{Q}_1}\left(\bar{x}^*(\delta)\right) \geq W_{Q_0}\left(-v^*(0) - \Delta x^*(0)\right) + W_{Q_1}\left(x^*(\delta)\right),$$

where $v^*(t), x^*(t)$ are defined as in Lemma 4.2.

*Proof.* By definition of support function, we can write

$$W_{\hat{Q}_1}\left(\bar{x}^*(\delta)\right) = \sup_{x(\delta) \in \hat{Q}_1} \langle x(\delta), \bar{x}^*(\delta)\rangle = \sup_{x(\delta) \in x(0) + \delta Q_1} \langle x(\delta), \bar{x}^*(\delta)\rangle$$

$$= \sup_{x(\delta) \in Q_0 + \delta Q_1} \langle x(\delta), \bar{x}^*(\delta)\rangle = W_{Q_0 + \delta Q}\left(\bar{x}^*(\delta)\right) = W_{Q_0}\left(\bar{x}^*(\delta)\right) + W_{\delta Q_1}\left(\bar{x}^*(\delta)\right)$$

$$= W_{Q_0}\left(\bar{x}^*(\delta)\right) + \delta W_{Q_1}\left(\bar{x}^*(\delta)\right) = W_{Q_0}\left(\bar{x}^*(\delta)\right) + W_{Q_1}\left(\delta \bar{x}^*(\delta)\right).$$

Therefore, taking into account that $\delta \bar{x}^*(0) = x^*(0)$, $\delta \bar{x}^*(\delta) = x^*(\delta)$, $\delta \bar{\mu}^*(0) = \mu^*(0)$, $\mu^*(0) = \delta v^*(0) + 2x^*(\delta)$ and recalling that the support function is positive homogeneous and subadditive, we derive the validity of the following relation

$$W_{Q_0}\left(\bar{x}^*(0) - \bar{\mu}^*(0)\right) + W_{\hat{Q}_1}\left(\bar{x}^*(\delta)\right) = W_{Q_0}\left(\bar{x}^*(0) - \bar{\mu}^*(0)\right) + W_{Q_1}\left(\delta \bar{x}^*(\delta)\right) + W_{Q_0}\left(\bar{x}^*(\delta)\right)$$

$$= W_{Q_0}\left(\frac{\delta \bar{x}^*(0) - \delta \bar{\mu}^*(0)}{\delta}\right) + W_{Q_1}\left(\delta \bar{x}^*(\delta)\right) + W_{Q_0}\left(\frac{\delta \bar{x}^*(\delta)}{\delta}\right)$$

$$= W_{Q_0}\left(\frac{x^*(0) - \mu^*(0)}{\delta}\right) + W_{Q_1}\left(x^*(\delta)\right) + W_{Q_0}\left(\frac{x^*(\delta)}{\delta}\right)$$

$$= W_{Q_0}\left(\frac{1}{\delta}\left[x^*(0) - \delta v^*(0) - 2x^*(\delta)\right]\right) + W_{Q_1}\left(x^*(\delta)\right) + W_{Q_0}\left(\frac{x^*(\delta)}{\delta}\right)$$

$$= W_{Q_0}\left(-v^*(0) - \Delta x^*(0) - \frac{x^*(\delta)}{\delta}\right) + W_{Q_1}\left(x^*(\delta)\right) + W_{Q_0}\left(\frac{x^*(\delta)}{\delta}\right)$$

$$\geq W_{Q_0}\left(-v^*(0)-\Delta x^*(0)\right)+W_{Q_1}\left(x^*(\delta)\right). \qquad \square$$

As consequence, using successively of Propositions 4.1-4.3 and Lemmas 4.2, 4.3 we can formulate the dual problem $(P_{DA}^*)$ for $(P_{DA})$ problem as follows

$$(P_{DA}^*) \quad \sup\left\{-\varphi^*\left(v^*(1-\delta)+\Delta_- x^*(1), -x^*(1)\right) + \sum_{t=0}^{1-2\delta} \delta M_{F(\cdot,t)}\left(\Delta^2 x^*(t)+\Delta_- v^*(t+\delta), v^*(t+\delta), x^*(t+2\delta)\right) \right.$$
$$\left. -W_{Q_0}\left(-v^*(0)-\Delta x^*(0)\right)-W_{Q_1}\left(x^*(\delta)\right)\right\},$$

where maximization is taken over the set $x^*(t), \mu^*(t), x^*(1), t=0,\ldots,1-\delta$. It should be noted that the main significance of the dual problem $(P_{DA}^*)$ is its expression in terms of difference derivatives.

## 5. The Dual Problem for Convex DFIs

In order to establish a dual problem to the main problem $(P_C)$, we use a limiting process in dual problem $(P_{DA}^*)$; by passing to the formal limit as $\delta \to 0$, the obtained maximization problem will be the dual problem to the previous continuous convex problem $(P_C)$. It is not hard ro see that here the second term in the problem $(P_{DA}^*)$ is an integral sums and as the mesh of the partition $\delta$ tends to zero, we have a definite integral of the integrant $M_{F(\cdot,t)}$ over time interval $[0,1]$:

$$(P_C^*) \quad \sup_{x^*(\cdot), v^*(\cdot)}\left\{-\varphi^*\left(v^*(1)+{x^*}'(1), -x^*(1)\right) + \int_0^1 \left[M_{F(\cdot,t)}\left({x^*}''(t)+{v^*}'(t), v^*(t), x^*(t)\right)\right]dt \right.$$
$$\left. -W_{Q_0}\left(-v^*(0)-{x^*}'(0)\right)-W_{Q_1}(x^*(0))\right\}.$$

Furthemore, we assume that $x^*(t)$, $t \in [0,1]$ is absolutely continuous function together with the first order derivatives and ${x^*}''(\cdot) \in L_1^n([0,1])$. Moreover, $v^*(\cdot)$ is absolutely continuous and ${v^*}'(\cdot) \in L_1^n([0,1])$.

In order that to prove the duality theorem we need to formulate the duality relation. We shall prove the duality relation is the Euler-Lagrange type adjoint inclusion. To this end, in the following theorem are formulated the sufficient conditions of optimality for the second order convex DFIs with convex initial point nonfunctional constraints $(P_C)$. These conditions are more precise since they involve useful forms of the Weierstrass-Pontryagin condition and second order Euler-Lagrange type adjoint inclusions. In the reviewed results this effort culminates in Theorem 5.1.

First, we formulate the reminded second order Euler-Lagrange type adjoint inclusion and transversality conditions for the problem $(P_C)$

$(a)$ $\left({x^*}''(t)+{v^*}'(t), v^*(t)\right) \in F^*\left(x^*(t); (\tilde{x}(t), \tilde{x}'(t), \tilde{x}''(t)), t\right)$, a.e. $t \in [0,1]$,

where

$(b)$ $\tilde{x}''(t) \in F_A\left(\tilde{x}(t), \tilde{x}'(t); x^*(t), t\right)$, a.e. $t \in [0,1]$.

The transversality conditions at the endpoints $t=0$ and $t=1$ consist of the following

$$(c) \quad \left(v^*(0) + x^{*\prime}(0), -x^*(0)\right) \in K_{Q_0}^*\left(\tilde{x}(0)\right) \times K_{Q_1}^*\left(\tilde{x}'(0)\right),$$

$$(d) \quad \left(v^*(1) + x^{*\prime}(1), -x^*(1)\right) \in \partial_{(x,y)} \varphi(\tilde{x}(1), \tilde{x}'(1)),$$

respectively. Now we are ready formulate the following theorem of optimality.

**Theorem 5.1** Suppose that $\varphi$ is a continuous and proper convex function, $F(\cdot, t)$ is a convex set-valued mapping and $Q_i (i = 0,1)$ are convex sets. Then for optimality of the feasible trajectory $\tilde{x}(t)$ in the problem $(P_C)$ it is sufficient that there exists a pair of absolutely continuous functions $\{x^*(t), v^*(t)\}$, $t \in [0,1]$ satisfying a.e. the second order Euler-Lagrange type differential inclusion $(a)$-$(b)$ and the transversality conditions $(c)$, $(d)$ at the initial point $t = 0$ and endpoint $t = 1$, respectively.

*Proof.* By the proof idea of Theorem 5.1[19] from $(a)$, $(b)$ we derive the following inequality

$$0 \geq \langle x'(1) - \tilde{x}'(1), x^*(1) \rangle - \langle x'(0) - \tilde{x}'(0), x^*(0) \rangle \\ - \langle v^*(1) + x^{*\prime}(1), x(1) - \tilde{x}(1) \rangle + \langle v^*(0) + x^{*\prime}(0), x(0) - \tilde{x}(0) \rangle. \tag{28}$$

Now, by definition of dual cones $K_{Q_0}^*\left(\tilde{x}(0)\right)$, $K_{Q_1}^*\left(\tilde{x}'(0)\right)$ from the transversality condition $(c)$ we deduce that

$$-\langle x'(0) - \tilde{x}'(0), x^*(0) \rangle + \langle v^*(0) + x^{*\prime}(0), x(0) - \tilde{x}(0) \rangle \geq 0, \forall x(0) \in Q_0; \forall x'(0) \in Q_1. \tag{29}$$

Thus, it follows from (28) and (29) that

$$0 \geq \langle x'(1) - \tilde{x}'(1), x^*(1) \rangle - \langle v^*(1) + x^{*\prime}(1), x(1) - \tilde{x}(1) \rangle. \tag{30}$$

Now, it is not hard to see that the transversality conditions $(d)$ at the endpoint $t = 1$, can be rewritten as follows

$$\varphi(x(1), x'(1)) - \varphi(\tilde{x}(1), \tilde{x}'(1)) \geq \langle v^*(1) + x^{*\prime}(1), x(1) - \tilde{x}(1) \rangle - \langle x^*(1), x'(1) - \tilde{x}'(1) \rangle. \tag{31}$$

Then, summing the inequalities (30), (31) for all feasible trajectories $x(\cdot)$, satisfying the initial conditions $x(0) \in Q_0$, $x'(0) \in Q_1$ we have the needed inequality:

$$\varphi(x(1), x'(1)) - \varphi(\tilde{x}(1), \tilde{x}'(1)) \geq 0 \quad \text{or} \quad \varphi(x(1), x'(1)) \geq \varphi(\tilde{x}(1), \tilde{x}'(1)). \qquad \square$$

We are now in a position to establish our duality relations between $(P_C)$ and $(P_C^*)$.

**Theorem 5.2** Suppose that the conditions of Theorem 4.1 are satisfied and $\tilde{x}(t)$ is an optimal solution of the primary problem $(P_C)$ with convex DFI. Then a pair of functions $\{\tilde{x}^*(\cdot), \tilde{v}^*(\cdot)\}$ is an optimal solution of the dual problem $(P_C^*)$ if and only if the conditions $(a)$–$(d)$ of Theorem 5.1 are satisfied. In addition, the optimal values in the primary $(P_C)$ and dual $(P_C^*)$ problems are equal.

□ Before all we prove that for all feasible solutions $x(\cdot)$ and dual variables $\{x^*(\cdot), v^*(\cdot)\}$ of the primary $(P_C)$ and dual $(P_C^*)$ problems, respectively, the inequality holds:

$$\varphi(x(1), x'(1)) \geq -\varphi^*\left(v^*(1) + x^{*\prime}(1), -x^*(1)\right) + \int_0^1 \left[M_{F(\cdot,t)}(x^{*\prime\prime}(t) + v^{*\prime}(t), v^*(t), x^*(t))\right] dt \tag{32}$$
$$-W_{Q_0}(-v^*(0) - x^{*\prime}(0)) - W_{Q_1}(x^*(0)).$$

To this end, by using the conjugate $\varphi^*$ and definition of Hamiltonian function we can write

$$\int_0^1 \left[M_{F(\cdot,t)}(x^{*\prime\prime}(t) + v^{*\prime}(t), v^*(t), x^*(t))\right] dt - W_{Q_0}\left(-v^*(0) - x^{*\prime}(0)\right) - W_{Q_1}(x^*(0))$$

$$-\varphi^*\left(v^*(1) + x^{*\prime}(1), -x^*(1)\right) \leq \int_0^1 \left[\langle x(t), x^{*\prime\prime}(t) + v^{*\prime}(t)\rangle + \langle x'(t), v^*(t)\rangle - \langle x''(t), x^*(t)\rangle\right] dt$$

$$+\varphi(x(1), x'(1)) - \langle x(1), v^*(1) + x^{*\prime}(1)\rangle - \langle x'(1), -x^*(1)\rangle + \langle v^*(0) + x^{*\prime}(0), x(0)\rangle$$

$$-\langle x^*(0), x'(0)\rangle = \int_0^1 \left[\langle x(t), x^{*\prime\prime}(t)\rangle - \langle x''(t), x^*(t)\rangle\right] dt + \int_0^1 d\langle x(t), v^*(t)\rangle$$

$$+\varphi(x(1), x'(1)) - \langle x(1), v^*(1) + x^{*\prime}(1)\rangle - \langle x'(1), -x^*(1)\rangle + \langle v^*(0) + x^{*\prime}(0), x(0)\rangle$$

$$-\langle x^*(0), x'(0)\rangle = \int_0^1 \left[\langle x(t), x^{*\prime\prime}(t)\rangle - \langle x''(t), x^*(t)\rangle\right] dt + \langle x(1), v^*(1)\rangle - \langle x(0), v^*(0)\rangle$$

$$+\varphi(x(1), x'(1)) - \langle x(1), v^*(1) + x^{*\prime}(1)\rangle - \langle x'(1), -x^*(1)\rangle + \langle v^*(0) + x^{*\prime}(0), x(0)\rangle \tag{33}$$

$$-\langle x^*(0), x'(0)\rangle = \int_0^1 \left[\langle x(t), x^{*\prime\prime}(t)\rangle - \langle x''(t), x^*(t)\rangle\right] dt + \varphi(x(1), x'(1)) - \langle x(1), x^{*\prime}(1)\rangle$$

$$-\langle x'(1), -x^*(1)\rangle + \langle x^{*\prime}(0), x(0)\rangle - \langle x^*(0), x'(0)\rangle.$$

Further, it is not hard to see that

$$\int_0^1 \left[\langle x^{*\prime\prime}(t), x(t)\rangle - \langle x''(t), x^*(t)\rangle\right] dt = \int_0^1 d\left[\langle x^{*\prime}(t), x(t)\rangle - \langle x'(t), x^*(t)\rangle\right]$$

$$= \langle x^{*\prime}(1), x(1)\rangle - \langle x'(1), x^*(1)\rangle - \langle x^{*\prime}(0), x(0)\rangle + \langle x'(0), x^*(0)\rangle. \tag{34}$$

Then the relationships (33) and (34) give us

$$\int_0^1 \left[M_{F(\cdot,t)}(x^{*\prime\prime}(t) + v^{*\prime}(t), v^*(t), x^*(t))\right] dt - W_{Q_0}\left(-v^*(0) - x^{*\prime}(0)\right) - W_{Q_1}(x^*(0))$$

$$-\varphi^*\left(v^*(1) + x^{*\prime}(1), -x^*(1)\right) \leq \langle x^{*\prime}(1), x(1)\rangle - \langle x'(1), x^*(1)\rangle - \langle x^{*\prime}(0), x(0)\rangle$$

$$+\langle x'(0), x^*(0)\rangle + \varphi(x(1), x'(1)) - \langle x(1), x^{*\prime}(1)\rangle - \langle x'(1), -x^*(1)\rangle$$

$$+\langle x^{*\prime}(0), x(0)\rangle - \langle x^*(0), x'(0)\rangle = \varphi(x(1), x'(1))$$

and this proves the inequality (32). Furthermore, suppose that a pair $\{\tilde{x}^*(\cdot), \tilde{v}^*(\cdot)\}$ satisfies the conditions (a)–(d) of Theorem 5.1. Then by definition of LAM the Euler-Lagrange type inclusion (a) and the condition (b) imply that

$$H_F\left(x(t), x'(t), \tilde{x}^*(t)\right) - H_F\left(\tilde{x}(t), \tilde{x}'(t), \tilde{x}^*(t)\right)$$
$$\leq \left\langle \tilde{x}^{*\prime\prime}(t) + \tilde{v}^{*\prime}(t), x(t) - \tilde{x}(t) \right\rangle + \left\langle \tilde{v}^*(t), x'(t) - \tilde{x}'(t) \right\rangle,$$

whence by the definition of function $M_F$ we deduce that

$$\left\langle \tilde{x}^{*\prime\prime}(t) + \tilde{v}^{*\prime}(t), \tilde{x}(t) \right\rangle + \left\langle \tilde{v}^*(t), \tilde{x}'(t) \right\rangle - H_F\left(\tilde{x}(t), \tilde{x}'(t), \tilde{x}^*(t)\right) \quad (35)$$
$$= M_F\left(\tilde{x}^{*\prime\prime}(t) + \tilde{v}^{*\prime}(t), \tilde{v}^*(t), \tilde{x}^*(t)\right).$$

On the other hand, by the transversality condition (c) we can write

$$-\left\langle \tilde{v}^*(0) + \tilde{x}^{*\prime}(0), \tilde{x}(0) \right\rangle = W_{Q_0}\left(-\tilde{v}^*(0) - \tilde{x}^{*\prime}(0)\right), \quad \left\langle \tilde{x}^*(0), \tilde{x}'(0) \right\rangle = W_{Q_1}\left(\tilde{x}^*(0)\right). \quad (36)$$

Finally, by Theorem 1.27 [15] the transversality condition (d) is equivalent to the relation

$$\varphi^*\left(\tilde{v}^*(1) + \tilde{x}^{*\prime}(1), -\tilde{x}^*(1)\right) = \left\langle \tilde{x}(1), \tilde{v}^*(1) + \tilde{x}^{*\prime}(1) \right\rangle + \left\langle \tilde{x}'(1), \tilde{x}^*(1) \right\rangle - \varphi\left(\tilde{x}(1), \tilde{x}'(1)\right). \quad (37)$$

Thus, taking into account the relationships (35)-37) in (33) the inequality sign is replaced by equality and for $\tilde{x}(\cdot)$ and $\{\tilde{x}^*(\cdot), \tilde{v}^*(\cdot)\}$ the equality of values of the primary and dual problems is ensured. Moreover, $\tilde{x}(\cdot)$ and $\{\tilde{x}^*(\cdot), \tilde{v}^*(\cdot)\}$ are satisfies the conditions (a)–(d) of Theorem 5.1 and the collection (a)–(d) is a dual relation for the primary $(P_C)$ and dual $(P_C^*)$ problems. □

**Example 5.1** Suppose we have the continuous case of Example 2.1 with semilinear second order DFI

$$\text{subject to} \quad \begin{aligned} &\text{infimum}\, \varphi(x(1), x'(1)), \\ &x''(t) = A_0 x(t) + A_1 x'(t) + Bu(t),\, u(\cdot) \in U, \\ &x(0) \in Q_0,\, x'(0) \in Q_1,\, t \in [0,1], \end{aligned} \quad (38)$$

where the matrices $A_0, A_1, B$, the function $\varphi$, and the sets $U, Q_i\,(i=0,1)$ are the same, as in Example 2.1. The problem is to find a controlling parameter $\tilde{u}(t) \in U$ such that the arc $\tilde{x}(t)$ corresponding to it minimizes $\varphi(x(1), x'(1))$. We introduce a set-valued mapping of the form $F(x, y) = A_0 x + A_1 y + BU$. Then using the formula for $M_{F(\cdot, t)}$ of Example.2.1, in view of the dual problem $(P_C^*)$ we can write

$$M_F\left(x^{*\prime\prime}(t) + v^{*\prime}(t), v^*(t), x^*(t)\right) = \begin{cases} -W_U(B^* x^*(t)), & \text{if } x^{*\prime\prime}(t) + v^{*\prime}(t) = A_0^* x^*(t), v^*(t) = A_1^* x^*(t), \\ -\infty, & \text{otherwise} \end{cases}$$

or, more compactly,

$$M_{F(\cdot,t)}\left(x^{*\prime\prime}(t) + v^{*\prime}(t), v^*(t), x^*(t)\right) = \begin{cases} -W_U(B^* x^*(t)), & \text{if } x^{*\prime\prime}(t) = A_0^* x^*(t) - A_1^* x^{*\prime}(t), \\ -\infty, & \text{otherwise.} \end{cases}$$

Then the dual problem of problem (38) is

$$\sup_{x^*(\cdot)}\left\{-\varphi^*\left(A_1^* x^*(1) + x^{*\prime}(1), -x^*(1)\right) - \int_0^1 W_U(B^* x^*(t))dt - W_{Q_0}\left(-A_1^* x^*(0) - x^{*\prime}(0)\right) - W_{Q_1}(x^*(0))\right\},$$

where $x^*(\cdot)$ is a solution of the adjoint Euler-Lagrange inclusion/equation $x^{*\prime\prime}(t) = A_0^* x^*(t)$

$-A_1^* x^{*\prime}(t)$. Consequently, maximization in this dual problem to primary problem (38) is realized over the set of solutions of the adjoint equation.

## 6. Conclusion

The paper under the "nondegeneracy" condition deals with the development of Mayer problem for second order evolution differential inclusions which are often used to describe various processes in science and engineering; the second-order discrete-approximate inclusions are investigated according to proposed discretization method; here we introduce a general model, that establishes a bridge between second order discrete and second order differential inclusions. First are derived necessary and sufficient optimality conditions in the form of Euler-Lagrange type inclusions and transversality conditions. Then we treat dual results according to the dual operations of addition and infimal convolution of convex functions. For construction of the duality problem skilfully computation of conjugate and support functions are required. It appears that the Euler-Lagrange type inclusions are duality relations for both prımary and dual problems and that the dual problem for discrete-approximate problem make a bridge between the dual problems of discrete and continuous problems. We believe that relying to the method described in this paper it can be obtained the similar duality results to optimal control problems with any higher order differential inclusions. In this way for computation of the conjugate function and support function of discrete-approximate problem a Pascal triangle with binomial coefficients, can be successfully used for any "higher order" calculations. These difficulties, of course, are connected with the existence of higher order difference derivatives in Mayer functional and discrete-approximate inclusions, respectively. Thus, the equivalence results for the conjugate functions and the Hamiltonian in the transition to the continuous problem are basic tools in the study of duality results; this approach plays a much more important role in construction of dual problems with second-order discrete and discrete-approximate inclusions. There has been a significant development in the study of duality theory to problems with first order differential/difference inclusions in recent years. As an open problem for further investigations, we mention the study of duality theory for an arbitrary higher-order differential inclusion. Besides, there can be no doubt that investigations of duality results to problems with second order differential inclusions can have great contribution to the modern development of the optimal control theory. Consequently, there arises a rather complicated problem with simultaneous determination of conjugacy of a Mayer functional depending of high order derivatives of searched functions. Thus, we can conclude that the proposed method is reliable for solving the various duality problems with higher order discrete and differential inclusions.

**Acknowledgments**. In advance, the author wishes to express his sincere thanks to Editor in-Chief, Prof. Irena Lasiecka of the Journal of Evolution Equations & Control Theory for consideration of this manuscript.

## References

[1] Artstein-Avidan and V. Milman, A characterization of the concept of duality, *Electronic Research Announcements,* **14** (2007), 42-59.

[2] D. Azzam-Laouir and F. Selamnia, On state-dependent sweeping process in Banach spaces, *Evol. Equ. Contr. Theory* (*EECT*), **7** (2018),[183-196.

[3] V. Barbu, I. Lasiecka, D. Tiba and  C. Varsan, Analysis and optimization of differential systems, IFIP TC7/WG7.2 International Worksing *Conference on Analysis and Optimization of Differential Systems, September 10-14, 2002, Constanta, Romania. IFIP*

*Conference Proceedings 249, Kluwer 2003,* ISBN 1-4020-7439-5.

[4] S.A. Belbas and S.M. Lenhart, Deterministic Optimal Control Problem with Final state Constraints, *Proceedings of the IEEE Conference On Decision and Control*. (1984), 526-527.

[5] A. Bressan, Differential inclusions and the control of forest fires, *J. Diff. Equ. (special volume in honor of A. Cellina and J. Yorke)*, **243** (2007), 179-207.

[6] G. Buttazzo and P.I Kogut, Weak optimal controls in coefficients for linear elliptic problems, *Revista Matematica Complutense*, **24** (2011), 83-94

[7] P. Cannarsa, A. Marigonda and K.T. Nguyen, Optimality conditions and regularity results for time optimal control problems with differential inclusions, *J. Math. Anal. Appl.* **427** (2015), 202-228.

[8] A. Dhara and A. Mehra, Conjugate duality for generalized convex optimization problems, J. Indust. Manag. Optim., **3** (2007), 415-427.

[9] M.D.Fajardol and J.Vidal, Necessary and sufficient conditions for strong Fenchel–Lagrange duality via a coupling conjugation scheme, *J. Optim. Theory Appl.* **176** (2018), 57-73.

[10] A.V. Fursikov, M. D. Gunzburger and L. S. Hou, Optimal boundary control for the evolutionary Navier-Stokes system: the three-dimensional case, *SIAM. J. Control Optim.* **43** (2005),2191-2232.

[11] A.D. Ioffe and V. Tikhomirov, Theory of extremal problems, "Nauka", Moscow, 1974 English transl., North-Holland, Amsterdam, 1978.

[12] N.C. Kourogenis, Strongly nonlinear second order differential inclusions with generalized boundary conditions, *J. Math. Anal. Appl.,* **287** (2003), 348-364.

[13] I. Lasiecka and N. Fourrier, Regularity and stability of a wave equation with strong damping and dynamic boundary conditions, *Evol. Equ. Contr. Theory (EECT),* **2** (2013), 631-667.

[14] P.J.Laurent, Approximation et optimisation, Herman, Paris, 1972.

[15] E.N. Mahmudov, Approximation and Optimization of Discrete and Differential Inclusions, Elsevier, Boston, USA , 2011.

[16] E.N. Mahmudov, On duality in problems of optimal control described by convex differential inclusions of Goursat-Darboux type, *J. Math. Anal. Appl.,* **307** (2005), 628-640.

[17] E.N. Mahmudov and M. E. Unal, Optimal control of discrete and differential inclusions with distributed parameters in the gradient form, *J. Dynam. Contr. Syst.,* 18 (2012), 83-101

[18] E.N. Mahmudov, Convex optimization of second order discrete and differential inclusions with inequality constraints. *J. Convex Anal.* **25** (2018), 1-26.

[19] E.N.Mahmudov, Optimization of Mayer problem with Sturm-Liouville type differential inclusions, *J. Optim.Theory Appl (JOTA),* **177** (2018), 345-375.

[20] E.N. Mahmudov, Optimization of fourth order Sturm-Liouville type differential inclusions with initial point constraints, *J. Industr. Manag. Optim. (JIMO),* (2018),13-35

[21] E.N. Mahmudov, Optimal control of second order delay-discrete and delay-differential inclusions with state constraints, *Evol. Equ. Contr. Theory (EECT),***7** (2018), 501-529.

[22] E.N.Mahmudov, Optimal control of higher order differential inclusions with functional constraints, *ESAIM: COCV*, doi: https://doi.org/10.1051/cocv/2019018.

[23] B.S. Mordukhovich and Tan H. Cao, Optimal control of a nonconvex perturbed sweeping process, *J. Differ. Equ.*, **266** (2019),1003-1050.

[24] N.S. Papageorgiou and Vicenţiu D. Rădulescu, Periodic solutions for time-dependent subdifferential evolution inclusions, *Evol. Equ. Contr. Theory (EECT),* **6** (2017), 277-297.

[25] R.T.Rockafellar and P.R. Wolenski, Convexity in Hamilton-Jacobi theory 1: Dynamics and duality, *SIAM J. Contr. Optim.*, **39** (2000), 1323-1350.

[26] T.I.Seidman, Compactness of a fixpoint set and optimal control, Appl. Anal., **88** (2009), 419-423.

[27] S. Sharma, A. Jayswal and S. Choudhury, Sufficiency and mixed type duality for multiobjective variational control problems involving α-V-univexity, *Evol. Equ. Contr. Theory* (*EECT*), **6**(2017), 93-109.

[28] Y. Zhou, V. Vijayakumar and R. Murugesu, Controllability for fractional evolution inclusions without compactness, *Evol. Equ. Contr. Theory* (*EECT*), **4** (2015), 507-524.

[29] Q.Zhang and G.Li, Nonlinear boundary value problems for second order differential inclusions, *Nonlin. Anal.: Theory, Methods Appl.*, **70** (2009), 3390-3406.

*E-mail address*:   elimhan22@yahoo.com